\documentclass[12pt, twoside, leqno]{article}
\usepackage{amsmath,amssymb,amsthm,latexsym}
\usepackage{enumerate}

\newtheorem{theorem}{Theorem}[section]

\newtheorem{lemma}[theorem]{Lemma}
\newtheorem{corollary}[theorem]{Corollary}

\date{}

\begin{document}

\title{Polynomial analogues of Ramanujan congruences for Han's hooklength formula}
\author{William J. Keith\\
Western Illinois University\footnote{This work was done during the author's postdoctorate at the University of Lisbon, Portugal.}\\
Email: WJ-Keith@wiu.edu}

\date{}

\maketitle

\renewcommand{\thefootnote}{}



\renewcommand{\thefootnote}{\arabic{footnote}}
\setcounter{footnote}{0}

\abstract{This article considers the eta power $\prod {(1-q^k)}^{b-1}$.  It is proved that the coefficients of $\frac{q^n}{n!}$ in this expression, as polynomials in $b$, exhibit equidistribution of the coefficients in the nonzero residue classes mod 5 when $n=5j+4$.  Other symmetries, as well as symmetries for other primes and prime powers, are proved, and some open questions are raised.}

\section{Introduction}

The Han/Nekrasov-Okounkov hook length formula, below, was discovered independently and by significantly different means, first by Nekrasov and Okounkov (\cite{NekOk}) and shortly thereafter by Guo-Niu Han (\cite{Han2008}).  The formula states that for any complex number $b$,

\begin{equation}\label{HNO}
\prod_{k=1}^{\infty} {(1-q^k)}^{b-1} = \sum_{\lambda \in {\cal{P}}} q^{\vert \lambda \vert} \prod_{h_{ij} \in \lambda} \left( 1 - \frac{b}{h_{ij}^2} \right) \, \text{ .}
\end{equation}

The formula relates powers of the partition generating function $\prod_{k=1}^{\infty} \frac{1}{1-q^k}$ or, inversely, the eta function $q^{\frac{1}{24}} \prod_{k=1}^{\infty} (1-q^k)$, to a sum of products over the hooklengths of all partitions.  The hooklengths $h_{ij}$ are a multiset of $n$ whole number values associated to any given partition of $n$; a detailed definition is delayed until a relevant question is proposed in the concluding section.  For any given $n$, the coefficient on $q^n$ is clearly a polynomial in the variable $b$.

Congruences for the numerical values of particular powers of the partition function have been studied (see \cite{AndMult},\cite{Atkin},\cite{CHL}, and \cite{Forbes}, among others; the first is a survey of related literature).  The object of this paper is to study the values of the coefficients of $q^n$ as polynomials in $b$.  After normalizing by $\frac{1}{n!}$ to make these coefficients integer polynomials, we find many pleasing results on the distribution and arrangement of the coefficients of these polynomials.

The most classical congruences for powers of the partition function are Ramanujan's congruences for the partition function itself, i.e. the case $b=0$.  Of particular interest to us is the first of these, that

$$p(5k+4) \equiv 0 \text{ mod } 5$$

\noindent for any integer $k$.  We show an analogue of this result for the hook length formula, namely,

\begin{theorem}\label{MainTheorem} For $n = 5k+4$, the integer polynomial $p_n (b)$ in the complex indeterminate $b$, defined by $\prod_{k=1}^{\infty} (1-q^k)^{b-1} = \sum_{n=0}^{\infty} \frac{q^n}{n!} p_n(b)$, has coefficients for which the following symmetries hold:

$\bullet$ The nonzero residues  equally populate the residue classes 1, 2, 3, and 4 mod 5.

$\bullet$ These residues appear in groups of four as a rotation of the list $(2,4,3,1)$.

$\bullet$ The coefficients of terms of degree at most $k$ are all 0 mod 5.  (There may be others.)
\end{theorem}

\vspace{0.1in}

\noindent \textbf{Remarks.} Concerning integrality, while it may not be immediately obvious that the polynomials are integral, this is easily shown by deriving a recurrence for the polynomials with the $q\frac{\partial}{\partial q} log$ technique of Herb Wilf's generatingfunctionology (\cite{WilfGF}), or referring to Corollary 2.3 of (\cite{Han08b}).\footnote{Because this recurrence is useful for calculating the polynomials discussed in this paper, we state it here for the benefit of readers.  With $p_0(b) = 1$, and $\sigma_1$ the divisor function, for $n\geq1$ we have $$p_n(b) = (n-1)!(b-1)\sum_{m=1}^n \frac{-\sigma_1(m) p_{n-m}(b)}{(n-m)!} \text{ .}$$  A little more work can derive an expression for the coefficients individually, giving the coefficient of $b^k$ in $p_n$ as an $(n-k) \times (n-k)$ (check this) determinant.}

\vspace{0.1in}

Equidistribution is only one of many beautiful symmetries exhibited by these coefficients.  Noting from the second clause of the theorem that the $4j$-th coefficient identifies coefficients through $4j+3$, we can can display every 4th coefficient starting from the $k+1$ position (the first $k$ being 0).  If we do this, leaving a space for zeroes for visual effect, we obtain the striking triangle below.  This is exactly Pascal's triangle, multiplied by 2 with alternating signs, reduced mod 5.

The structure of this paper is as follows.  In the next section we prove Theorem 1.1.  In Section 3 we show additional symmetries, which are corollaries of the proof of Theorem 1 and well-known facts on binomial coefficients.  We then consider other prime and prime power arithmetic progressions.  In Section 4 we discuss a few open questions on the combinatorics of these polynomials, which will hopefully motivate future work in this study.

For exposition suggestions and catching certain typographical errors the author cordially thanks the referee.

\vspace{0.1in}

{\tt \begin{tabular}{ll}
n=4: &  \{2\} \\
n=9: &  \{2,3\} \\
n=14: &  \{2,1,2\} \\
n=19: &  \{2,4,1,3\} \\
n=24: &  \{2,2,2,2,2\} \\
n=29: &  \{2,\phantom{-},\phantom{-},\phantom{-},\phantom{-},3\} \\
n=34: &  \{2,3,\phantom{-},\phantom{-},\phantom{-},3,2\} \\
n=39: &  \{2,1,2,\phantom{-},\phantom{-},3,4,3\} \\
n=44: &  \{2,4,1,3,\phantom{-},3,1,4,2\} \\
n=49: &  \{2,2,2,2,2,3,3,3,3,3\} \\
n=54: &  \{2,\phantom{-},\phantom{-},\phantom{-},\phantom{-},1,\phantom{-},\phantom{-},\phantom{-},\phantom{-},2\} \\
n=59: &  \{2,3,\phantom{-},\phantom{-},\phantom{-},1,4,\phantom{-},\phantom{-},\phantom{-},2,3\} \\
n=64: &  \{2,1,2,\phantom{-},\phantom{-},1,3,1,\phantom{-},\phantom{-},2,1,2\} \\
n=69: &  \{2,4,1,3,\phantom{-},1,2,3,4,\phantom{-},2,4,1,3\} \\
n=74: &  \{2,2,2,2,2,1,1,1,1,1,2,2,2,2,2\} \\
n=79: &  \{2,\phantom{-},\phantom{-},\phantom{-},\phantom{-},4,\phantom{-},\phantom{-},\phantom{-},\phantom{-},1,\phantom{-},\phantom{-},\phantom{-},\phantom{-},3\} \\
n=84: &  \{2,3,\phantom{-},\phantom{-},\phantom{-},4,1,\phantom{-},\phantom{-},\phantom{-},1,4,\phantom{-},\phantom{-},\phantom{-},3,2\} \\
n=89: &  \{2,1,2,\phantom{-},\phantom{-},4,2,4,\phantom{-},\phantom{-},1,3,1,\phantom{-},\phantom{-},3,4,3\} \\
n=94: &  \{2,4,1,3,\phantom{-},4,3,2,1,\phantom{-},1,2,3,4,\phantom{-},3,1,4,2\} \\
n=99: &  \{2,2,2,2,2,4,4,4,4,4,1,1,1,1,1,3,3,3,3,3\} \\
n=104: &  \{2,\phantom{-},\phantom{-},\phantom{-},\phantom{-},2,\phantom{-},\phantom{-},\phantom{-},\phantom{-},2,\phantom{-},\phantom{-},\phantom{-},\phantom{-},2,\phantom{-},\phantom{-},\phantom{-},\phantom{-},2\} \\
n=109: &  \{2,3,\phantom{-},\phantom{-},\phantom{-},2,3,\phantom{-},\phantom{-},\phantom{-},2,3,\phantom{-},\phantom{-},\phantom{-},2,3,\phantom{-},\phantom{-},\phantom{-},2,3\} \\
n=114: &  \{2,1,2,\phantom{-},\phantom{-},2,1,2,\phantom{-},\phantom{-},2,1,2,\phantom{-},\phantom{-},2,1,2,\phantom{-},\phantom{-},2,1,2\} \\
n=119: &  \{2,4,1,3,\phantom{-},2,4,1,3,\phantom{-},2,4,1,3,\phantom{-},2,4,1,3,\phantom{-},2,4,1,3\} \\
n=124: &  \{2,2,2,2,2,2,2,2,2,2,2,2,2,2,2,2,2,2,2,2,2,2,2,2,2\} \\
n=129: &  \{2,\phantom{-},\phantom{-},\phantom{-},\phantom{-},\phantom{-},\phantom{-},\phantom{-},\phantom{-},\phantom{-},\phantom{-},\phantom{-},\phantom{-},\phantom{-},\phantom{-},\phantom{-},\phantom{-},\phantom{-},\phantom{-},\phantom{-},\phantom{-},\phantom{-},\phantom{-},\phantom{-},\phantom{-},3\} \\
n=134: &  \{2,3,\phantom{-},\phantom{-},\phantom{-},\phantom{-},\phantom{-},\phantom{-},\phantom{-},\phantom{-},\phantom{-},\phantom{-},\phantom{-},\phantom{-},\phantom{-},\phantom{-},\phantom{-},\phantom{-},\phantom{-},\phantom{-},\phantom{-},\phantom{-},\phantom{-},\phantom{-},\phantom{-},3,2\} \end{tabular}} \hfill

\rm 

\section{Proof of the main theorem}

We begin by expanding the product out using the generalized binomial theorem.  Recall that for any value $x$, including the indeterminate in $\mathbb{C}[x]$, we can define the generalized binomial coefficient with a whole number $k$ as $\binom{x}{k} = \frac{x(x-1)(x-2)\cdots(x-k+1)}{k!}$.

We use the notation $\mathbf{e} \vdash n$ to mean that $\mathbf{e}$ is a partition of $n$, and write partitions in the frequency notation $\mathbf{e}=1^{e_1}2^{e_2}3^{e_3}\dots$ to denote the partition in which 1 occurs $e_1$ times, 2 occurs $e_2$ times, etc.  It is understood that the $e_i$ are nonnegative integers and that only finitely many of the $e_i$ are nonzero.

Expanding with the generalized binomial theorem, we first obtain

\begin{multline*}\prod_{j=1}^{\infty} {(1-q^j)}^{b-1} = \prod_{j=1}^{\infty} \frac{1}{(1-q^j)^{1-b}} \\
 = \prod_{j=1}^{\infty} \sum_{k=0}^{\infty} \binom{1-b+k-1}{k} {(q^j)}^k  = \prod_{j=1}^{\infty} \sum_{k=0}^{\infty} \binom{k-b}{k} {(q^j)}^k \\
 = \prod_{j=1}^{\infty} \sum_{k=0}^{\infty} \frac{(k-b)(k-b-1)\cdots(k-b-k+1)}{k!} {(q^j)}^k \\
 = \prod_{j=1}^{\infty} \sum_{k=0}^{\infty} \frac{(1-b)(2-b)\cdots(k-b)}{k!} {(q^j)}^k \\
 = \sum_{n=0}^{\infty} q^n \sum_{{{\mathbf{e} \vdash n} \atop {\mathbf{e}} = 1^{e_1}2^{e_2}\dots}} \prod_{j=1}^{\infty} \frac{(1-b)(2-b)\cdots(e_j-b)}{e_j!} \\
 = \sum_{n=0}^{\infty} \frac{q^n}{n!} \sum_{{{\mathbf{e} \vdash n} \atop {\mathbf{e}} = 1^{e_1}2^{e_2}\dots}}  \frac{n!}{e_1!e_2!\dots} \prod_{j=1}^{\infty} (1-b)(2-b)\dots(e_j-b) \, \text{.}
 \end{multline*}
 
Note that $ \frac{n!}{e_1!e_2!\dots}$ is not a multinomial coefficient since $e_1 + e_2 + \dots \neq n$ except in the trivial case $e_1 = n$.  The $k=0$ term in the sums, an empty product, is 1.
 
 Considering the various powers of $b$ in the final polynomials, one polynomial for each $e_i$, we see that this expression is equal to
 
$$\sum_{n=0}^{\infty} \frac{q^n}{n!} \sum_{{{\mathbf{e} \vdash n} \atop {\mathbf{e}} = 1^{e_1}2^{e_2}\dots}} \frac{n!}{e_1!e_2!\dots} \prod_{j=1}^{\infty} \sum_{t=0}^{e_j} (-b)^{t} \sum_{{S=\{s_1,\dots,s_{e_j-t}\}} \atop {S \subseteq \{1,\dots,e_j\}}} s_1 s_2 \cdots s_{e_j-t} \text{ .}$$

The products on the far right are every possible selection of $e_j-t$ distinct elements from the set $\{1,2,\dots,e_j\}=: [e_j]$.  When we bring the factor of $\frac{1}{e_j!}$ back into this product, we exchange this for a sum over their complements, all possible denominators consisting of products of $t$ distinct elements from $\{1,2,\dots,e_j\}$.  

Multiplying these together, we get
 
$$\prod_{k=1}^{\infty} (1-q^k)^{b-1} = \sum_{n=0}^{\infty} \frac{q^n}{n!} \!\!\! \sum_{{{\mathbf{e} \vdash n} \atop {\mathbf{e}} = 1^{e_1}2^{e_2}\dots}} \!\!\! n!  \left[ \sum_{t=0}^n (-b)^t \left(\sum \frac{1}{s_1 \dots s_t} \right) \right]$$

\noindent where the last sum runs over sets of $t$ distinct elements chosen from the multiset $M_{\mathbf{e}} = \{1,2,\dots,e_1,1,2,\dots,e_2,\dots\}$.

Thus, the coefficient of $b^t$ in $p_n(b)$ is

$$(-1)^t n! \sum_{{{\mathbf{e} \vdash n} \atop {\mathbf{e}} = 1^{e_1}2^{e_2}\dots}} \left(\sum_{{S \subseteq M_{\mathbf{e}}} \atop {S = \{ s_1,\dots,s_t\}}}  \frac{1}{s_1 \dots s_t} \right) \, \text{ .}$$

Our task is now to determine the residue class of this integer mod 5 when $n=5k+4$.

In any given term, if the power of 5 that divides $n!$ is not fully canceled by elements of the product $s_1 \dots s_t$, that term will contribute 0 to the residue class of the sum mod 5.  It is possible for this to occur if $e_1 \geq 5k$, and among the $s_i$ are $5, 10,15,\dots,5k$ chosen from the first part of the multiset, $\{1,2,\dots,e_1,(\dots)\}$.  If $e_1 < 5k$, since all other parts are of size at least 2, it is clear that the deficit in available entries can never be fully made up by elements chosen from $\{1,\dots,e_j\}$ with $j>1$.  For example, if $e_2 \geq 5$, then $e_1 < 5k-5$, etc.

Hence the only partitions $\mathbf{e}$ that can possibly contribute to the sum mod 5 are:

$$1^{5k+4} \quad , \quad 1^{5k+2}2^1 \quad , \quad 1^{5k+1}3^1 \quad , \quad 1^{5k} 4^1 \quad , \quad 1^{5k} 2^2 \, \text{ .}$$

Among the contributions from these partitions, only those terms in which $k$ of the elements are assigned to $5, 10, \dots, 5k$ can possibly contribute to the residue of the sum.

Thus, if $t < k$, no terms contribute a nonzero value.  If $t=k$, then exactly the one term $\frac{(5k+4)!}{5\cdot10\cdot\dots\cdot5k}$, from $S=\{5,10,\dots,5k\}$, contributes.  This happens once in each of the five possible partitions, so the total is 0 mod 5.  This proves the last clause of Theorem 1.

Suppose now that $t=k+m$, $m>0$, and that $k$ of the elements of $\{s_1,\dots,s_t\}$ are $5,10,\dots,5k$.  There remain the contributions constructed from choosing $m$ values from the remaining places in the five possible partitions.  Since we are only concerned with the residue of the result mod 5, we need only classify the multiset of available choices in each partition by their residues mod 5.

If we reduce the elements of the multiset $M_{\mathbf{e}}$ mod 5 for each of the five possible $\mathbf{e}$, we obtain the following five multisets:

\begin{itemize}
\item $1^{5k+0}4^1$: $\{\overbrace{1,2,3,4,\phantom{5},}^{k \text{ repetitions}}\dots,1\}$
\item $1^{5k+0}2^2$: $\{\overbrace{1,2,3,4,\phantom{5},}^{k \text{ repetitions}}\dots,1,2\}$
\item $1^{5k+1}3^1$: $\{\overbrace{1,2,3,4,\phantom{5},}^{k \text{ repetitions}}\dots,1,1\}$
\item $1^{5k+2}2^1$: $\{\overbrace{1,2,3,4,\phantom{5},}^{k \text{ repetitions}}\dots,1,2,1\}$
\item  $1^{5k+4}$: $\{\overbrace{1,2,3,4,\phantom{5},}^{k \text{ repetitions}}\dots,1,2,3,4\}$
\end{itemize}

Let the part of each possible $M_{\mathbf{e}}$ consisting of $k$ repetitions of $\{1,2,3,4\}$ be $C$.  Since $C$ is common to all five partitions, if $S \subseteq C \subset M_{\mathbf{e}}$ for some $\mathbf{e}$ then the same choice may be made in any of the other four partitions.  The summand $$\frac{(5k+4)!}{5\cdot10\cdot\dots\cdot5k} \frac{1}{s_{k+1} \dots s_t}$$ \noindent thus appears five times in the coefficient on $b^t$.  The sum of these five terms is then 0 mod 5, and we may consider the remaining possibilities for $S$.  Our strategy is now to form more collections of terms which vanish mod 5, and determine the remaining total.

Each $M_{\mathbf{e}}$ has a 1 immediately following $C$.  If we choose $m-1$ elements of $S$ from $C$ and the first following 1, our choice of $S$ can be found in each $M_{\mathbf{e}}$, and the terms $$\frac{(5k+4)!}{5\cdot10\cdot\dots\cdot5k}\frac{1}{s_{k+1} \dots s_{t-1}} \frac{1}{1}$$ can be grouped to vanish mod 5.

We again group all those terms where $m-1$ elements of $S$ are the same choices made in $C$, and the additional 2 is chosen -- there are three of these -- or the second 1 of the partition $1^{5k+1}3^1$ is chosen.  These give us

\begin{multline*}
3 \cdot \frac{(5k+4)!}{5\cdot10\cdot\dots\cdot5k}\frac{1}{s_{k+1} \dots s_{t-1}} \frac{1}{2} + \frac{(5k+4)!}{5\cdot10\cdot\dots\cdot5k}\frac{1}{s_{k+1} \dots s_{t-1}} \frac{1}{1} \\
\equiv \frac{(5k+4)!}{5\cdot10\cdot\dots\cdot5k}\frac{1}{s_{k+1} \dots s_{t-1}} \left(5 \cdot 2^{-1} \right)  \, \text{ mod } \, 5.
\end{multline*}

These thus contribute 0 residue mod 5.

We again group all those terms where $m-1$ of the elements of $S$ are in $C$, and the second 1 of $1^{5k+2}2^1$ is chosen, or the 4 of $1^{5k+4}$.  We obtain 

\begin{multline*}
\frac{(5k+4)!}{5\cdot10\cdot\dots\cdot5k}\frac{1}{s_{k+1} \dots s_{t-1}} \frac{1}{1} + \frac{(5k+4)!}{5\cdot10\cdot\dots\cdot5k}\frac{1}{s_{k+1} \dots s_{t-1}} \frac{1}{-1} \\
\equiv \frac{(5k+4)!}{5\cdot10\cdot\dots\cdot5k}\frac{1}{s_{k+1} \dots s_{t-1}} \left(0\right)  \, \text{ mod } \, 5.
\end{multline*}

Thus, \emph{among those terms in which $m-1$ of the elements of $S$ are in $C$}, we are left with those in which the 3 of the $1^{5k+4}$ is chosen, multiplying the contribution of these elements by $3^{-1} \equiv 2$ mod 5.

What is the actual value of that contribution?  If $n=5k+4$, then there are $k$ 1s, $k$ 2s, et cetera.  All possible $S$ with $m-1$ choices in $C$ thus yield $$\frac{(-1)^t(5k+4)!}{5\cdot10\cdot\dots\cdot5k} \sum_{r_1+r_2+r_3+r_4=m-1} \binom{k}{r_1}\binom{k}{r_2}\binom{k}{r_3}\binom{k}{r_4} {\frac{1}{1}}^{r_1}{\frac{1}{2}}^{r_2}{\frac{1}{3}}^{r_3}{\frac{1}{4}}^{r_4}$$

\noindent where the $r_i$ index those choice of $S$ in which $r_i$ elements of $S$ are among the residues $i$ in $C$.  This value will then be multiplied by $3^{-1}$ mod 5 in the final total.

Now suppose $m-2$ elements of $S$ are in $C$.  We shall be briefer, since the argument is similar.

We can form the following large group that contributes 0 mod 5: any choice of 1 and 2 from the same set (there are 3 of these), or the 2 and the second 1 of $1^{5k+2}2^1$, the pairs of 1s from $1^{5k+1}3^1$ and $1^{5k+2}2^1$, and of $1^{5k+4}$ the 1 and 3, the 1 and 4, the 2 and 4, and the 3 and 4.  Ignoring the prior factor, the sum of the inverses is

\begin{multline*}3 \times 1^{-1}2^{-1} + 2^{-1}1^{-1} + 2 \times 1^{-1}1^{-1} + 1^{-1}3^{-1} + 1^{-1}4^{-1} + 2^{-1}4^{-1} + 3^{-1}4^{-1}\\
\equiv 9+3+2+2-1-3-2 = 10 \equiv 0 \, \text{ mod } \, 5.
\end{multline*}

\emph{Among those terms in which $m-2$ of the choices are made in the first $k$ sets of 4}, we are left with those in which the 2 and 3 of the $1^{5k+4}$ is chosen, multiplying the contribution of these elements by 1 mod 5.

If $m-3$ choices are made in the first $k$ sets of 4, then we have 3 choices to make from the remaining elements.  The 4 possible choices from $1^{5k+4}$ -- $\{1,2,3\}$, $\{1,2,4\}$, $\{1,3,4\}$, and $\{2,3,4\}$ -- group to add 0 mod 5.  We are left with the choice of 1, 2, and 1 from $1^{5k+2}2^1$, mutliplying the contribution of these choices by 3.

If $m-4$ choices are made in the first $k$ sets of 4, the only remaining possibility is all of 1,2,3, and 4 from $1^{5k+4}$, multiplying these contributions by 4.

For any given $t=k+m$, the sum is thus a linear combination of the residues mod 5 of the various ways to choose $m-1$ through $m-4$ elements of $S$ from $C$.  We have therefore reduced the question of determining the residue class of the coefficient of $b^t$ in $p_n(b)$ to calculating the residue class mod 5 of the following sum:

\begin{multline*}\frac{(-1)^t(5k+4)!}{5\cdot10\cdot\dots\cdot5k} \left( \sum_{r_1+r_2+r_3+r_4=m-1} \frac{1}{3}\binom{k}{r_1}\binom{k}{r_2}\binom{k}{r_3}\binom{k}{r_4} {\frac{1}{1}}^{r_1}{\frac{1}{2}}^{r_2}{\frac{1}{3}}^{r_3}{\frac{1}{4}}^{r_4} \right. \\
+ \sum_{r_1+r_2+r_3+r_4=m-2} \frac{1}{6}\binom{k}{r_1}\binom{k}{r_2}\binom{k}{r_3}\binom{k}{r_4} {\frac{1}{1}}^{r_1}{\frac{1}{2}}^{r_2}{\frac{1}{3}}^{r_3}{\frac{1}{4}}^{r_4} \\
+ \sum_{r_1+r_2+r_3+r_4=m-3} \frac{1}{2}\binom{k}{r_1}\binom{k}{r_2}\binom{k}{r_3}\binom{k}{r_4} {\frac{1}{1}}^{r_1}{\frac{1}{2}}^{r_2}{\frac{1}{3}}^{r_3}{\frac{1}{4}}^{r_4} \\
\left. + \sum_{r_1+r_2+r_3+r_4=m-4} \frac{1}{24}\binom{k}{r_1}\binom{k}{r_2}\binom{k}{r_3}\binom{k}{r_4} {\frac{1}{1}}^{r_1}{\frac{1}{2}}^{r_2}{\frac{1}{3}}^{r_3}{\frac{1}{4}}^{r_4} \right) \, \text{ .}
\end{multline*}

Taking inverses mod 5, the coefficient of $b^t$ in $p_n(b)$, where $t=k+m$ and $n=5k+4$,  becomes

$$(-1)^{m+1} \left( \sum_{{r_1+r_2+r_3+r_4=m-c}\atop{c=1,2,3,4}} a_c\binom{k}{r_1}\binom{k}{r_2}\binom{k}{r_3}\binom{k}{r_4} {1}^{r_1}{3}^{r_2}{2}^{r_3}{4}^{r_4} \right)$$

\noindent where $a_1 = 2$, $a_2 = 1$, $a_3 = 3$, and $a_4 = 4$.

We now wish to evaluate this sum, whose symmetries over the range of $m$ range yield the behaviors of Theorem 1 and the additional symmetries visible in the triangle.  It will be useful to evaluate the sum for a general prime $p$, as follows:

\begin{lemma}\label{SumLemma} For $p$ a prime, \begin{multline*}\sum_{r_1+\dots + r_{p-1}=(p-1)s+c} \binom{k}{r_1} \cdot \dots \cdot \binom{k}{r_{p-1}} {1}^{r_1}{2}^{r_2}\dots{(p-1)}^{r_{p-1}} \\ \equiv \left\{ \begin{matrix} (-1)^s \binom{k}{s} \, \text{ mod } p & c=0 \\ 0 \, \text{ mod } p & \text{ otherwise.} \end{matrix} \right. \end{multline*}
\end{lemma}

\noindent \textbf{Proof.} (Note that we assign bases to exponents of their own index, for convenience.  The sum is symmetric under this exchange.) 

The left-hand side is the coefficient of $q^{(p-1)s+c}$ in the product \begin{multline*}(1+q)^k \dots (1+(p-1)q)^k = (1 + (\sum_{i=1}^{p-1}i)q + (\sum_{{i,j=1} \atop {i<j}}^{p-1}ij)q^2 + \dots + (p-1)!q^{p-1})^k \\
\equiv (1 - q^{p-1})^k \, \text{ mod } \, p \, \text{ .}
\end{multline*}

The middle terms in the sum vanish (observe the effect of permutation by a nontrivial multiplication in $\mathbb{Z}_p$) and since $(p-1)! \equiv -1$ mod $p$, we obtain the last line.  The coefficient of $q^{(p-1)s}$ in the latter expression is exactly $(-1)^s \binom{k}{s}$, and all other coefficients are 0.  \hfill $\Box$

We now complete the proof of the first two clauses of theorem \ref{MainTheorem}.  The coefficients of $b^t$ are sums of multiples of four consecutive terms of the type evaluated in Lemma \ref{SumLemma}, only one of which may be nonzero mod 5.  The corresponding coefficients of degree $k+4s+1$, $k+4s+2$, $k+4s+3$ and $k+4s+4$ are this value times 2, 1, 3, and 4, respectively, and an alternating power of $-1$ with starting parity determined by $s$.

Thus the rotation is the underlying value multiplied by either $(2,-1,3,-4)=(2,4,3,1)$, or by $(-2,1,-3,4)=(3,1,2,4)$.  Both are rotations of the list $(2,4,3,1)$.  If an underlying nonzero value is $j$ mod 5, then the list $(2j,4j,3j,1j)$ is a rotation of $(2,4,3,1)$, as claimed.  Finally, when the underlying sum is zero mod 5, all 4 terms are 0 mod 5.  

This means that those coefficients on $b^t$ which are nonzero rotate through the nonzero residue classes mod 5 in groups of 4, leading to equidistribution of the pattern claimed in Theorem \ref{MainTheorem}. $\Box$

\section{Additional Symmetries}

We have now evaluated the coefficient of $b^{k+1+4m}$ in $p_{5k+4} (b)$ to be congruent to $2 \cdot (-1)^m \binom{k}{m}$ mod 5.  This tells us that the entries of the triangle possess all the symmetries of Pascal's triangle, when multiplied by 2 mod 5 and an alternating sign.

For example, we know that for any prime $p$ the power of $p$ in the prime factorization of the binomial coefficient $\binom{k}{m}$ is equal to the number of carries in the addition $m + (k-m) = k$ in base-$p$ arithmetic.  Thus

\begin{corollary}  The coefficient of $b^{k+1+4m}$ in $p_{5k+4} (b)$, $0 \leq m \leq k$, is nonzero if and only if $m_i \leq k_i$ for each digit in the 5-ary expansions $(m)_5 = m_0m_1\dots$ and $(k)_5 = k_0k_1\dots$ \, .
\end{corollary}

Modulo any prime, Pascal's triangle possesses self-similarity:

\begin{corollary} The triangle of Section 1 is self-similar: the apex digits of the triangles at any level are given by the entries of the triangle one level down, and other entries are those of the fundamental triangle times the apex digit.
\end{corollary}

Similar corollaries of any of the symmetries of Pascal's triangle reduced mod 5 can be used.  One must take the multiple $2 (-1)^m$ into account, of course; the coefficients of $p_n(b)$ are not fixed under reversal, for instance.

\subsection{Other Primes}

It is not the case that equidistribution holds for the nonzero residues classes modulo other primes, whether in the $-1$ arithmetic progression or in those for which Ramanujan-like congruences hold for the partition function.  For example, $p_6(b) = 7920 - 18144b + 14674 b^2 - 5205 b^3 + 805 b^4 - 51 b^5 + b^6$.  However, mod 7, hese coefficients reduce to $(3,0,2,3,0,5,1)$ and likewise $p_5(b) \equiv (0,6,4,2,0,6)$.

This occurs because the $a_c$ are not as neatly distributed for other primes as they are for modulus 5.  However, equidistribution mod $p$ arises for a different reason in the arithmetic progression $-p-1$ mod $p^2$.  This is because the binomial coefficients described by Lemma \ref{SumLemma} themselves rotate through the nonzero residues mod $p$.

\begin{theorem} For $p$ prime, $j \geq 1$, if the number of partitions of $p-1$ is not congruent to 1 mod $p$, the coefficients of $p_{p^2-p-1}(b)$ equinumerously populate the nonzero residue classes mod $p$ for all $j$, and if it is, the populations are still equinumerous for $p_{p^3-p^2-p-1}$ mod $p$.\end{theorem}

In order to prove the theorem for all primes, we need to go through the proof of Theorem 1 in greater generality.

The expansion makes no reference to the value of $p$ until we begin determining the residue class of the coefficient $$(-1)^t n! \sum_{{{\mathbf{e} \vdash n} \atop {\mathbf{e}} = 1^{e_1}2^{e_2}\dots}} \left(\sum_{{S \subseteq M_{\mathbf{e}}} \atop {S = \{ s_1,\dots,s_t\}}}  \frac{1}{s_1 \dots s_t} \right) \, \text{ .}$$

The same logic holds to show that the number of 1s $e_1$ must be sufficient to cancel all instances of the prime $p$ in $n!$.  Thus for any given prime $p$ and residue $r$, if we consider the arithmetic progress $n=pk+r$, only the partitions with $e_1 \geq pk$ contribute, and there are a finite set of these equal in number to the partitions of $r$.

We can run the same grouping arguments to give a series of coefficients $a_c$, $0 \leq c \leq p-1$.  The coefficient of $b^t$ in $p_n(b)$, where $n = p(pj+p-2)+(p-1)$ and $t = pj+p-2+m$, $m \geq 0$, is congruent mod $p$ to

\begin{multline*}(-1)^{m+1} \left( \sum_{{r_1+\dots+r_{p-1}=m-c} \atop {c=0,1,\dots,p-1}} a_c \binom{k}{r_1}\dots\binom{k}{r_{p-1}} {1}^{r_1}{2}^{r_2}\dots{(p-1)}^{r_{p-1}} \right) \, \text{ .}
\end{multline*}

This is a linear combination of terms of the form addressed by Lemma \ref{SumLemma}, so only every $p-1$st term is nonzero.  It was the case in the proof of Theorem \ref{MainTheorem} that $a_0=0$ and so there was no overlap at the ends of this interval, where a term might take contributions from both $a_0$ and $a_4$, but this might not be the case for other primes.  Mod $p$ the sequence of coefficients of $p_n(b)$, with the initial segment of 0s of length $pj+p-2$, is given by

\begin{multline*} (0,\dots,0,-a_0 \binom{pj+p-2}{0}, \\ a_1\binom{pj+p-2}{0},\dots,a_{p-2}\binom{pj+p-2}{0}, \hfill \\ \hfill -a_{p-1} \binom{pj+p-2}{0} + a_0 \binom{pj+p-2}{1}, \\ -a_1\binom{pj+p-2}{1}, \dots , -a_{p-2}\binom{pj+p-2}{1}, \hfill \\ \hfill a_{p-1} \binom{pj+p-2}{1} - a_0 \binom{pj+p-2}{2}, \\ a_1\binom{pj+p-2}{2},\dots, (-1)^{pj+2p-3} a_{p-2}\binom{pj+p-2}{pj+p-2}, \hfill \\ \hfill (-1)^{pj+2p-2} a_{p-1} \binom{pj+p-2}{pj+p-2}+ (-1)^{pj+p} a_0\binom{pj+p-2}{pj+p-1}) \text{ .}
\end{multline*}

There are two types of subsequences here: those multiplied by $a_c$ for $1 \leq c \leq p-2$, which are of the form $\{ (-1)^{c+1} a_c (-1)^s \binom{pj+p-2}{s}\}$, $0 \leq s \leq pj+p-2$, and the overlapping sums for $a_0$ and $a_{p-1}$.  Let us begin by examining the former.

It is a straightforward calculation to show that

\begin{lemma}\label{BinomLemma} For $p$ a prime, $0 \leq s \leq pj+p-2$, $s = gp+h$, $0 \leq h < p$, $$(-1)^s \binom{pj+p-2}{s} \equiv (h+1) \left( (-1)^g \binom{j}{g} \right) \, \text{ mod } p \, \text{ .}$$ \end{lemma}

This shows that the sequences $\{ (-1)^{c+1} a_c (-1)^s \binom{pj+p-2}{s}\}$ consist of segments of $p$ elements that are either all 0, or are permutations of $\{1,2,\dots,p-1\}$ followed by a 0.  (The last 0 is in the ``place'' $s = pj+p-1$.)  Thus the nonzero elements equally populate the residue classes mod $p$

Using Lemma \ref{BinomLemma}, the sequence of overlapping terms at the ends of the intervals, $\{-a_0 (-1)^s \binom{pj+p-2}{s} - a_{p-1} (-1)^{s-1} \binom{pj+p-2}{s-1} \}$, where $0 \leq s \leq pj+p-1$, $s = gp+h$ with $0 \leq h< p$, reduces to 

\begin{multline*} \{ -(a_0(h+1)\left( (-1)^g \binom{j}{g} \right) + a_{p-1} (h) \left( (-1)^{\lfloor \frac{s-1}{p} \rfloor} \binom{j}{\lfloor \frac{s-1}{p} \rfloor} \right) ) \} \equiv \\
\{ (-1)^0 \binom{j}{0} (-a_0), (-1)^0 \binom{j}{0} (-2a_0-1a_{p-1}), (-1)^0 \binom{j}{0} (-3a_0-2a_{p-1}), \phantom{\}} \\
\dots,(-1)^0 \binom{j}{0} (-p a_0 -(p-1)a_{p-1}), a_0 \binom{j}{1}, (-1)^1 \binom{j}{1} (-2a_0-a_{p-1}), \\
\phantom{\{} \dots , (-1)^j\binom{j}{j}(-p a_0 - (p-1) a_{p-1} ) \} \, \text{ .}
\end{multline*} 

This is the set of values $(-1)^x\binom{j}{x} (-a_0 + y(-a_0 - a_{p-1}))$ where $0 \leq x \leq j$ and $0 \leq y \leq p-1$.  If $a_0 + a_{p-1} \not\equiv$ 0 mod $p$, this will equally populate all nonzero residue classes mod $p$ each time $y$ ranges over its values.

Now we note that $a_{p-1}$ is always -1 mod $p$, because the only contributing partition must be $1^{pk+p-1}$ and only the terms in which all of 1 through $p-1$ beyond $C$ are chosen.  Further, $a_0$ is the number of partitions of $p-1$, because this is exactly the number of possible contributing partitions.  Hence $a_0 + a_{p-1} \not\equiv 0$ mod $p$ whenever the number of partitions of $p-1$ is not congruent to 1 mod $p$.

Even when $a_0 + a_{p-1} \equiv 0$ mod $p$, by Lemma \ref{BinomLemma} the sequence will still equally populate the nonzero residue classes mod $p$ as $x$ ranges over its values when $j \equiv -2$ mod $p$, giving us the arithmetic progression $-p^2-p-1$ mod $p^3$. \hfill $\Box$

\vspace{0.1in}

\noindent \textbf{Remarks:} The first time the exceptional case happens is at $p=71$.  The size of the populations in the arithmetic progression $-1-p-p^2-p^3-p^4-\dots$ mod $p^q$ for any $q$ can be seen, from the argument above, to be divisible by $(p-1)^{q-2}$, or $(p-1)^{q-3}$ in the exceptional case.

\section{Open Questions}

This sequence of polynomials is related to a core object in enumerative combinatorics, so questions regarding its symmetries should be of wide interest, and many still await investigation.

One of the most interesting questions is whether we can recover the results of the earlier work mentioned on congruences of powers of the partition function, and perhaps extend them, by means as elementary as possible.  If we assume knowledge of the fact the the number of partitions of $5k+4$ is itself divisible by 5, it seems likely that the theorems and lemmas in this paper, combined with elementary properties of binomial coefficients, could yield useful theorems, such as the fact that $\prod (1-q^j)^{b-1}$ is divisible by 0 in the arithmetic progression $n=5k+4$ as long as $ b \not\equiv 3$ mod 5.  Similar results might be obtainable for other prime power progressions.  Exploration of these ideas is intended as the immediate followup to this paper.

Are there any other prime progressions $pk+(p-1)$ where equidistribution occurs?  If so, how can we find them efficiently?  If not, how could this be proved?

What can we say about progressions with composite moduli other than prime powers?

Is it possible to describe an evocative combinatorial object that the coefficients themselves count?  Could this description be useful in proofs regarding their properties, or properties of related objects such as multipartitions?

Finally, consider the Han/Nekrasov-Okounkov formula which was the original object of this investigation.  If we equate that expression for the $p_n(b)$ with the one obtained in this paper, we get

$$\sum_{\lambda \in {\cal{P}}} \prod_{h_{ij} \in \lambda} \left( 1 - \frac{b}{h_{ij}^2} \right) = \sum_{{{\lambda \vdash n} \atop {\lambda} = 1^{e_1}2^{e_2}\dots}} \!\!\! \left[ \sum_{t=0}^n (-b)^t \left(\sum_{{S \subseteq M_{\lambda}} \atop {S = \{ s_1,\dots,s_t\}}} \frac{1}{s_1 \dots s_t} \right) \right] .$$

This identity definitely does \emph{not} refine to the individual-partition level, yet the $M_\lambda$ are merely a subset of the hooklengths $h_{ij}$.  We finally find it necessary to define the hooklengths of a partition more rigorously -- if we write $\lambda = (\lambda_1, \lambda_2,\dots,\lambda_\ell)$, with the parts in nonincreasing order, then the hooklengths $h_{ij}$ are the multiset of values $h_{ij} = \lambda_i -j + \#\{\lambda_a \vert a \geq i \text{ and } \lambda_a \geq j\} $.  The elements of $M_\lambda$ are just the hooklength $h_{i{\lambda_i}}$, i.e. the "top strip" of hooklengths.

Thus one might ask, is there any refinement of the above identity?  For instance, the contribution of the unique 1-row partition $1^n$ is the same on both sides.  Could something be said about partitions with a given number of rows?

\end{document}